\documentclass[11pt]{amsart}
\usepackage{breqn,color}
\usepackage{amsmath,amsthm,amssymb,amsfonts,mathtools,bm,bbm,array,float,mathdots,dsfont,mathrsfs,latexsym}
\usepackage{caption,subcaption,multirow,longtable,lscape,wrapfig,tikz,tikz-cd,comment,pgffor,booktabs}
\usepackage[colorlinks=true, pdfstartview=FitV,linkcolor=blue,citecolor=blue,urlcolor=blue]{hyperref}
\usepackage{url}
\usepackage{fullpage}
\usepackage{xcolor}
\usepackage[colorinlistoftodos]{todonotes}
\usepackage{verbatim}
\usepackage{numprint}
\npthousandsep{,}  
\allowdisplaybreaks

\usepackage[normalem]{ulem}
\usepackage[style=alphabetic, backend=biber, sorting=nyt, doi=false, isbn=false,url=false, hyperref,backref,backrefstyle=none, maxnames=20, maxalphanames=20]{biblatex}

\usepackage{afterpage}

\definecolor{revblue}{RGB}{0,0,130}
\newcommand{\rev}[1]{{\color{revblue}#1}}

\newtheorem {theorem}    {Theorem}[section]

\theoremstyle{definition}

\newtheorem{convention}[theorem]{Conventions}

\newenvironment{red}{\relax\color{red}}{\hspace*{.5ex}\relax}
\newenvironment{blue}{\relax\color{revblue}}{\hspace*{.5ex}\relax}

\newcommand{\ber}{\begin{red}}
\newcommand{\er}{\end{red}}
\newcommand{\beb}{\begin{blue}}
\newcommand{\eb}{\end{blue}}
\newcommand{\seteq}{\mathbin{:=}}

\numberwithin{equation}{section}
\numberwithin{figure}{section}
\numberwithin{table}{section}
\allowdisplaybreaks

\begin{document}

\title{Murmurations, Mestre--Nagao sums, and Convolutional Neural Networks for elliptic curves}

\author[J. Bieri]{Joanna Bieri}
\address{University of Redlands, Redlands, CA 92373, USA}
\email{\href{mailto:joanna_bieri@redlands.edu}{joanna\_bieri@redlands.edu}}

\author[E. Costa]{Edgar Costa}
\address{
  Department of Mathematics,
  Massachusetts Institute of Technology,
  Cambridge,
  MA 02139,
  USA
}
\email{\href{mailto:edgarc@mit.edu}{edgarc@mit.edu}}
\urladdr{\url{https://edgarcosta.org}}

\author[A. Deines]{Alyson Deines}
\address{Center for Communications Research, La Jolla, USA}
\email{\href{mailto:aly.deines@gmail.com}{aly.deines@gmail.com}}

\author[K.-H. Lee]{Kyu-Hwan Lee}
\address{Department of Mathematics, University of Connecticut, Storrs, CT 06269, USA  \hfill \break \indent Korea Institute for Advanced Study, Seoul 02455, Republic of Korea}
\email{\href{mailto:khlee@math.uconn.edu}{khlee@math.uconn.edu}}

\author[D. Lowry-Duda]{David Lowry-Duda}
\address{ICERM, Providence, RI, 02903, USA
\hfill \break \indent
Department of Mathematics, Brown University, Providence, RI, 02912, USA}
\email{\href{mailto:david@lowryduda.com}{david@lowryduda.com}}
\urladdr{\url{https://davidlowryduda.com}}

\author[T. Oliver]{Thomas Oliver}
\address{University of Westminster, London, UK}
\email{\href{mailto:T.Oliver@westminster.ac.uk}{T.Oliver@westminster.ac.uk}}

\author[Y. Qi]{Yidi Qi}
\address{Department of Physics, Northeastern University, Boston, MA, USA \hfill \break \indent NSF Institute for Artificial Intelligence and Fundamental Interactions, Cambridge, MA, USA \hfill \break \indent School of Engineering and Applied Sciences, Harvard University, Cambridge, MA, US}
\email{\href{mailto:yidi_qi@seas.harvard.edu}{yidi_qi@seas.harvard.edu}}

\author[T. Veenstra]{Tamara Veenstra}
\address{Center for Communications Research, La Jolla, USA}
\email{\href{mailto:tamarabveenstra@gmail.com}{tamarabveenstra@gmail.com}}

\date{\today}

\begin{abstract}
We apply one-dimensional convolutional neural networks to  the Frobenius traces of elliptic curves over $\mathbb{Q}$ and evaluate and interpret their predictive capacity.
In keeping with similar experiments by Kazalicki--Vlah, Bujanovi\'{c}--Kazalicki--Novak, and Pozdnyakov, we observe high accuracy predictions for the analytic rank across a range of conductors. 
We interpret the prediction using saliency curves and explore the interesting interplay between saliency, murmurations and Mestre--Nagao sums.
\end{abstract}

\maketitle

\section{Introduction}
The vanishing order of a rational $L$-function at its central point is a fundamental arithmetic invariant.
For example, this vanishing order for the $L$-function associated with an abelian variety is expected to be equal to the rank of the finitely generated abelian group of rational points.
Inspired partly by this connection to a mysterious arithmetic quantity, a recent theme in the literature has been to machine learn the vanishing order of a rational $L$-function from a finite list of (normalized) Frobenius traces indexed by primes \cite{HLOc,  KV22, HLOP,Poz, many2}.  
This body of work fits into a wider literature of machine learning arithmetic quantities including class numbers \cite{HLOb, AHLOS}, Galois groups \cite{LL25}, Tate--Shafarevich groups \cite{BBFHS} and Fricke signs \cite{many}.
All told, it is well documented that standard machine learning algorithms can achieve high accuracy in predicting arithmetic quantities.

A more interesting, and more challenging, question is what to make of the apparent success of machine learning in this context. For example, can highly accurate arithmetic predictions be attributed to meaningful arithmetic structures, and can these structures be uncovered by interpreting how the predictions are made?
A preliminary step in this direction led to the discovery of so-called {\em murmurations} \cite{HLOP}.
Although originating from experimental work, murmurations are a genuinely number-theoretic phenomenon that can be studied rigorously by methods independent of machine learning \cite{Zub23, LOPdirichlet, BLLDSHZ, LDtrace}.
At the time of discovery, it was noted heuristically that the murmuration patterns associated with a dataset of elliptic curves appear correlated with the weightings in the first principal component of that dataset \cite{HLOP}.
This suggests, perhaps counter-intuitively, that when elliptic curves with conductor in a fixed interval are separated by rank using principal component analysis, the Frobenius traces at certain primes contribute more significantly than others.
Subsequently, Pozdnyakov observed murmurations by applying interpretability techniques to convolutional neural networks trained to predict analytic rank, more precisely, by analyzing their learned convolutional filters \cite{Poz}. 
At the same time, Pozdnyakov observed that Mestre--Nagao sums also contribute to rank predictions, which is perhaps unsurprising since they are expected to converge to a value closely connected to the rank, and they have historically played a significant role in the discovery of elliptic curves with high rank \cite{EK20}.
In contrast to the oscillating feature importance suggested by murmurations, Mestre--Nagao sums weight primes up to a fixed bound, with feature importance decaying on the order of 
$\log(p)/p$.
On the other hand, in keeping with the observation that certain primes may be more pertinent than others to rank prediction, it has been observed that Mestre--Nagao sums up to smaller bounds can sometimes yield more accurate predictions than those up to larger ones, and that the numerical discrepancy can be attributed to murmurations \cite{BKN}. 

In this paper, we apply a different interpretability technique to a CNN trained to predict analytic rank on a larger dataset of elliptic curves.
We reproduce earlier observations that both murmurations and Mestre–Nagao sums influence rank predictions for elliptic curves made by convolutional neural networks (CNNs). We then extend this analysis by systematically investigating their relative impact. Whereas Pozdnyakov analyzes binned averages of $a_p$ (which amplifies the impact of murmurations), we are able to examine the individual values.
This constitutes our first main contribution.
Our primary tool is the study of {\em saliency} curves for the CNNs, which are analogous to the principal component weightings considered in \cite{HLOP}, but contrast with the convolutional filters analyzed in \cite{Poz}.
The evolution of these saliency curves reveals how the balance between murmuration effects and Mestre–Nagao sums changes over the course of training. 
In early epochs, the saliency curves exhibit clear murmuration-like patterns, whereas in later epochs the importance of primes decays in a manner resembling the Mestre–Nagao weightings.

We follow this up by examining how the comparative size of the primes and conductors affects this balance.
In particular, murmuration patterns are observable as a function of $p/N$, where $N$ is the conductor \cite{sarnakLetter}, and this is also visible in the rank predictions.
For example, using only primes $p \leq 10^4$, murmuration-like patterns are visible only in the saliency curves for conductor up to $10^4$, and, for these conductors, the murmurations play a substantially larger role in rank prediction.
However, when we increase the range to $p \leq B$ for each conductor interval $[A, B]$ and retrain the models, we find that murmuration patterns emerge in the saliency curves at every conductor range we consider.
On the other hand, we also observe that the models trained on larger primes are slightly less accurate than those trained on $p\le 10^4$ alone.
More generally, the contribution of murmuration phenomenon to rank prediction depends on several factors, including the training epoch, the prime range, and the rank itself. 
For example, the murmuration-like patterns are clearest in the saliency curves associated with ranks $0$ and $1$.

We conclude this introduction with an overview of the following sections.
In Section~\ref{ss.SL_ECQ}, we describe our datasets and the architecture of our CNN, report prediction accuracy, and investigate how the accuracy depends on the number of primes.
In Section~\ref{s.saliency}, we present the saliency curves, and explore the interplay among predicted rank, murmurations, Mestre--Nagao sums, and conductor ranges.

\subsection*{Acknowledgments} We thank Alexey Pozdnyakov for helpful discussions, and the anonymous referee for their careful reading of the manuscript, which led to significant improvements. We also thank the Harvard University Center of Mathematical Sciences and Applications and, in particular, the organizers and attendees of the Mathematics and Machine Learning Program from 2024 and the Reunion from 2025. This collaboration grew out of these programs.
Costa was supported by Simons Foundation grant SFI-MPS-Infrastructure-00008651.
Qi was supported by the NSF grant PHY-2019786 (the NSF AI Institute for Artificial Intelligence and Fundamental Interactions).
We also acknowledge Claude Code Opus 5 for assistance with coding, scripting, and proof reading. 

\section{Data, architecture, and accuracy}\label{ss.SL_ECQ}

In this section, we seek to predict the rank of elliptic curves using one-dimensional CNNs. We specify our dataset in Section~\ref{ss.ecd} and our CNN architecture in Section~\ref{s.cnna}.
We achieve accuracy comparable to that in \cite{KV22} and \cite{Poz}, and record how the accuracy varies with the number of Frobenius traces used to represent the data points.

\subsection{Data specification}\label{ss.ecd}
Let $E$ be an elliptic curve over $\mathbb{Q}$. 
If $p$ is a good prime for $E$, then we define $a_p(E) \seteq p+1-\#E(\mathbb{F}_p)$. 
If $p$ is a bad prime for $E$, then  $a_p(E)\in\{-1,0,1\}$ is defined according to the reduction type for $E\bmod p$ \cite{silverman-advanced}.
The \emph{conductor} of $E$ is a positive integer whose prime divisors are precisely the bad primes, with exponents determined by the reduction type at each such prime \cite{Silverman:2009}. 
Roughly speaking, the conductor measures the arithmetic complexity of the elliptic curve.
According to the Hasse bound, for a prime $p$, we may normalise $a_p$ as follows:
\[\widetilde{a}_p(E)=\frac{a_p(E)}{2\sqrt{p}}\in[-1,1].\]
In this section, we will use a dataset called \texttt{XECQ}, which contains $\numprint{1536784}$ elliptic curves over $\mathbb{Q}$ with conductor no greater than $\numprint{400000}$
(\texttt{XECQ} is a subset of the larger set, available at~\cite{ecq6ap}).
Each elliptic curve in \texttt{XECQ} 
has (analytic) rank in the set $\{0,1,2,3,4\}$.
We represent each datapoint $E$ in \texttt{XECQ} by a vector of the form: 
\begin{equation}\label{eq.vectors}
v_E=(\widetilde{a}_p(E))_{p \le b}\in[-1,1]^{\pi(b)},
\end{equation}
where $b\in\mathbb{Z}_{>0}$ and $\pi(b)$ denotes the number of primes $\leq b$.
For example, if $b=\numprint{10000}$, then $\pi(b)=\numprint{1229}$.
We allow $b$ to vary with the conductor.

\begin{figure}[htbp]
    \centering
    \includegraphics[width=0.7\textwidth]{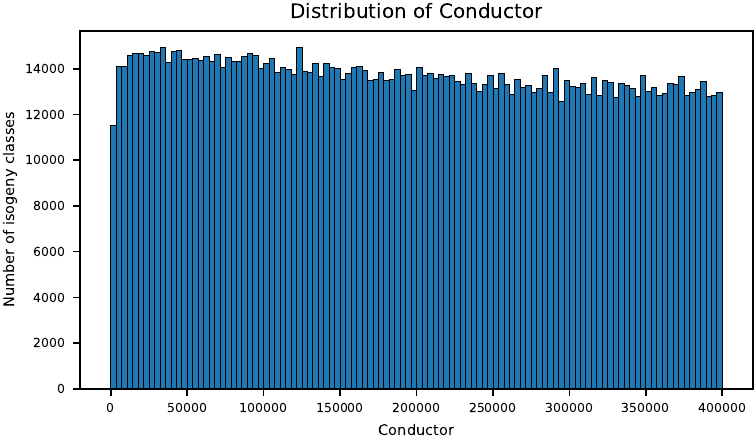}
    \captionof{figure}{Distribution of conductors in \texttt{XECQ}.}
    \label{fig:conductor_distribution} 
\end{figure}

The distribution of conductors within \texttt{XECQ} is shown in Figure \ref{fig:conductor_distribution}.
For an interval $[A,B]\subset\mathbb{R}_{>0}$, we denote by $\texttt{XECQ[A,B]}$ the subset of $\texttt{XECQ}$ containing only curves with conductor in $[A,B]$. 
In Table~\ref{tab:rank_counts}, we display four choices of intervals $[A,B]$ that yield approximately the same number of data points in $\texttt{XECQ[A,B]}$, and record the distribution of analytic ranks for each.
For every curve in these $\texttt{XECQ[A,B]}$ we compute $a_p$ for all primes $p\leq B$.

\begin{table}[H]
    \centering
    \begin{tabular}{|c|c|c|c|c|c|}
        \hline
        $[A,B]$ & $\#\texttt{XECQ[A,B]}$ & rank $0$ & rank $1$ & rank $2$ & rank $3$ \\ \hline
        $[0, 10000]$       & 
        $\numprint{36837}$ &
        $\numprint{15766}$ & $\numprint{19102}$ & $\numprint{1968}$ & $1$ \\ \hline
        $[100000, 110000]$ & 
        $\numprint{39667}$ &
        $\numprint{14575}$ & $\numprint{19869}$ & $\numprint{5082}$ & $141$ \\ \hline
        $[200000, 210000]$ & 
        $\numprint{38905}$ &
        $\numprint{13711}$ & $\numprint{19692}$ & $\numprint{5322}$ & $180$ \\ \hline
        $[300000, 310000]$ & 
        $\numprint{36742}$ &
        $\numprint{12793}$ & $\numprint{18374}$ & $\numprint{5345}$ & $230$ \\ \hline
    \end{tabular}
    \captionof{table}{Number of curves of each analytic rank in each \texttt{XECQ[A,B]}. }
    \label{tab:rank_counts}
\end{table}

\subsection{Architecture and accuracy}\label{s.cnna}
The code used for the experiments described in this paper is available at \cite{code-repo}.
For each of the four intervals $[A,B]$ listed in Table~\ref{tab:rank_counts}, we split $\texttt{XECQ[A,B]}$ into training, validation and test sets using an $80{:}10{:}10$ ratio and train a one-dimensional convolutional neural network (1D~CNN) to predict the rank.
We train for at most $300$ epochs, stopping early when the validation accuracy has not improved for $30$ epochs, and we report the test accuracy at the epoch of maximal validation accuracy.

The architecture of our neural network is summarised in Table~\ref{tab:cnn_architecture}.
In particular, the network consists of three 1D convolutional layers, with $16$, $32$ and $64$ channels respectively, each using kernel size $K = 3$ and padding $P = 1$. 
Each convolutional layer is followed by a $\mathrm{ReLU}$ activation function, and a max-pooling layer with $K = 2$ and $P = 1$. 
After a dropout layer with $D = 0.5$, the convolutional blocks are followed by two fully connected layers with $128$ neurons each, and a final output layer with width $5$ (that is, the number of possible vanishing orders in \texttt{XECQ}).
The network is trained with cross-entropy loss using the Adam optimizer, with batch size $\numprint{3000}$ and learning rate $0.001$.
Each epoch consists of ten steps based on the chosen mini-batch size.

\begin{table}[htbp]
    \centering
    \begin{tabular}{llc}
        \toprule
        Layer (Type) & Output Shape & Param \# \\
        \midrule
        Conv1D & $[1, 16, \pi(b)]$ & 64 \\
        MaxPool1D & $[1, 16, 615]$ & -- \\
        Conv1D & $[1, 32, 615]$ & 1,568 \\
        MaxPool1D & $[1, 32, 308]$ & -- \\
        Conv1D & $[1, 64, 308]$ & 6,208 \\
        MaxPool1D & $[1, 64, 155]$ & -- \\
        Flatten \& Dropout & $[1, 9920]$ & -- \\
        Linear & $[1, 128]$ & 1,269,888 \\
        Linear & $[1, 128]$ & 16,512 \\
        Linear & $[1, 4]$ & 516 \\
        \midrule
        \multicolumn{2}{l}{Total Trainable Parameters} & 1,294,756 \\
        \bottomrule
    \end{tabular}
    \caption{Architecture summary of the 1D CNN.}
    \label{tab:cnn_architecture}
\end{table}

In summary, our model is a shallow, fixed-depth architecture that progressively increases the number of channels, applies max-pooling, incorporates fully connected layers, and uses explicit dropout regularization. In contrast, the CNN of~\cite{KV22} is a modular, fully convolutional architecture
that increases the number of channels only once, employs strided convolutions for
downsampling in dedicated layers, and preserves feature dimensions elsewhere.
For larger-scale experiments, Pozdnyakov adopts a similar architecture to that
of~\cite{KV22}, as described in~\cite[Section~4.2]{Poz}.

In Table~\ref{tab:accuracy_XECQ}, we report the test accuracy for each conductor interval $[A,B]$, both for the fixed bound $b=10^4$ and for conductor-dependent bound $b=B$ (for $[A,B]=[0,10000]$, the two bounds coincide).
As in~\cite{KV22} and~\cite{Poz}, we observe high predictive accuracy.
For the conductor-dependent bound, the number of input features is $\pi(B)$, that is, $\numprint{1229}$, $\numprint{10453}$, $\numprint{18807}$ and $\numprint{26800}$ for the four intervals respectively.
Complete confusion matrices are available with the code \cite{code-repo}.
Notably, extending the prime range to align with conductor magnitude does not improve accuracy: the models trained on $p\le 10^4$ are more accurate in every interval. 
Experiments using different random seeds give similar results to those in Table~\ref{tab:accuracy_XECQ}.
Indeed, the $b=10^4$ test accuracies agree within $0.2$ percentage points and the $b=B$ accuracies agree within $2$ percentage points, and in all cases the $b=10^4$ models are more accurate than the $b=B$ models.

\begin{figure}[htbp]
    \centering
    \includegraphics[width=0.7\textwidth]{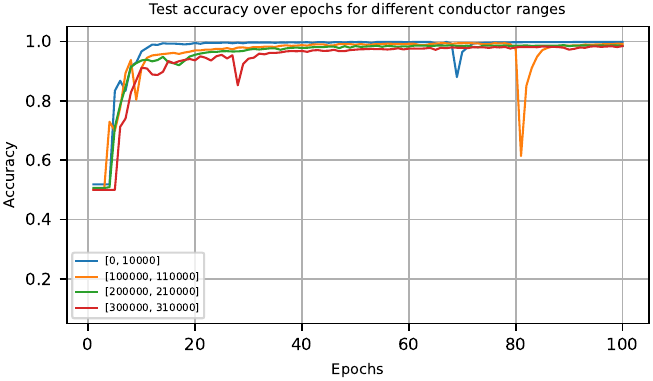}
    \captionof{figure}{Learning the rank in \texttt{XECQ}[A,B]: test accuracy against epoch. Each colour corresponds to one of the $4$ conductor intervals specified.}
    \label{fig:loss_XECQ}
\end{figure}

Figure~\ref{fig:loss_XECQ} displays the prediction accuracy of each CNN as a function of the number of training epochs.
Two of the curves in Figure~\ref{fig:loss_XECQ} exhibit sharp transient drops in accuracy. 
These drops are optimization artifacts, unrelated to the data.
In later training stages, CNNs can sit in a local minimum with only small gradients. This can prevent the model from leaving the local minimum and getting better overall accuracy. The Adam optimizer makes occasional large steps to try to escape these local minima, corresponding to the transient drops. After a few epochs the model recovers. The nature of which epoch and runs experience these drops, as well as the size of the drops, depends on random seeding.

\begin{figure}[htbp]
    \centering
    \begin{subfigure}[b]{0.45\textwidth}
        \centering
    \includegraphics[width=\textwidth]{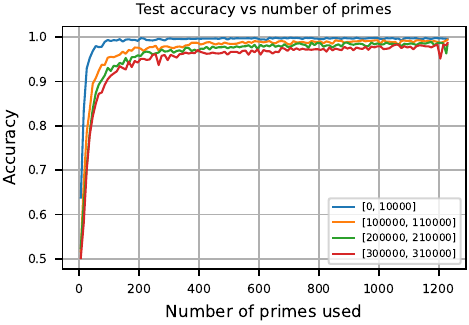}
    \end{subfigure}
    \hfill
    \begin{subfigure}[b]{0.45\textwidth}
        \centering
        \includegraphics[width=\textwidth]{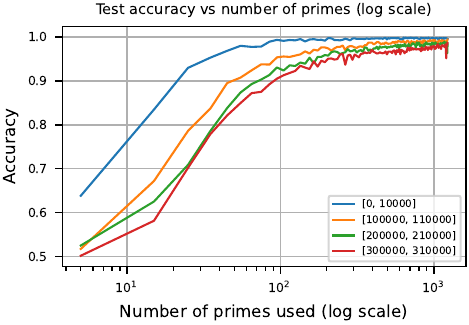}
    \end{subfigure}
    \caption{Learning the rank in \texttt{XECQ}[A,B]: (Left) test accuracy against number of primes used in training for different choices of $[A,B]$. (Right) test accuracy against the logarithm of the number of primes used in training  for different choices of $[A,B]$. }
    \label{fig:accuracy_vs_n_ap}
\end{figure}

It is generally expected that curves with larger conductor are more arithmetically complex and therefore may require more information to achieve the same level of accuracy in rank prediction.
With this in mind, the left (resp.\ right) panel of Figure~\ref{fig:accuracy_vs_n_ap}
plots the prediction accuracy as a function of $\pi(b)$ (resp.\ $\log \pi(b)$),
that is, the number of features (resp.\ its logarithm).

\begin{table}
    \begin{tabular}{|c|c|c|}
        \hline
        $[A,B]$ & test accuracy, $b = 10^4$ & test accuracy, $b = B$ \\ \hline
        $[0, 10000]$ &  $99.81\%$ & --- \\ \hline
        $[100000, 110000]$  & $99.27\%$ & $98.54\%$ \\ \hline
        $[200000, 210000]$ & $98.74\%$ & $97.53\%$\\ \hline
        $[300000, 310000]$ & $98.34\%$ & $96.63\%$\\ \hline
    \end{tabular}
    \captionof{table}{Test accuracy at the epoch of maximal validation accuracy for the 1D CNN described in Section~\ref{s.cnna} when applied to \texttt{XECQ}[A,B], for the fixed bound $b=10^4$ and the conductor-dependent $b=B$.}
    \label{tab:accuracy_XECQ}
\end{table}

\section{Saliency curves}\label{s.saliency}
To interpret the predictions of the CNN from Section~\ref{s.cnna}, we analyze the feature importance using saliency scores.
In \cite{Poz}, two key influences on analytic rank prediction are identified: murmurations and Mestre--Nagao sums.
In this section, we will observe that the relative contributions of these factors vary depending on the predicted rank, the training stage, and the range of primes shown to the network relative to the conductor.
When the models are trained on the fixed prime range $p \le 10^4$, the murmuration contribution appears to diminish with growing conductor. 
On the other hand, when they are trained on primes up to conductor bound, murmuration patterns emerge at every conductor range we consider, and they persist longer in training.

\subsection{Definitions}
After each epoch, the CNN from Section~\ref{s.cnna} determines a function
\begin{equation}\label{eq.CNN}
[-1,1]^{\pi(b)}\rightarrow\{0,1,2,3\},\end{equation}
which assigns to each input vector $(x_q)_{q \le b}$, with $x_q \in [-1,1]$ for primes $q \le b$, a predicted class (that is, predicted rank).
Although the function in~\eqref{eq.CNN} depends on the epoch, we will suppress this from subsequent notation.
We note that, although there are rank $4$ curves with conductor less than our largest conductor bound \numprint{310000}, there are no rank 4 curves in the sets $\texttt{XECQ}[A,B]$ used in this paper.

The function in equation~\eqref{eq.CNN} factors through a map
\[
S : [-1,1]^{\pi(b)} \longrightarrow \mathbb{R}^4,
\]
\[
S((x_q)_{q \le b}) = (S_0(\mathbf x),S_1(\mathbf x),S_2(\mathbf x),S_3(\mathbf x)),
\]
where, for each class $v \in \{0,1,2,3\}$, the quantity $S_v(\mathbf x)$ denotes the score assigned to the input $\mathbf x \seteq (x_q)_{q \le b}$.
The probability that a given input $\mathbf x \in [-1,1]^{\pi(b)}$ belongs to class $v$ is then defined by the softmax function
\[
    P_v (\mathbf x) = \frac{e^{S_v(\mathbf x)}}{\sum_{i=0}^{4} e^{S_i(\mathbf x)}},
\]
where $v \in \{0,1,2,3,4\}$.
The predicted class of the input $\mathbf x$ is given by
\[
c (\mathbf x) = \operatorname*{argmax}_{v \in \{0,1,2,3,4\}} P_v(\mathbf x) .
\]
Consider a class $v$ and an input $\mathbf x =(x_q)_{q \le b}\in [-1,1]^{\pi(b)}$. For  a prime $p \le b$, define the {\em saliency score} $w_p^v (\mathbf x)$ to be  the derivative of $S_v$ with respect to $x_p$ at $\mathbf x \in [-1,1]^{\pi(b)}$: 
\begin{equation}\label{eq.wp-1}
     w_p^v(\mathbf x) = \frac{\partial S_v (\mathbf x)}{\partial x_p}, \qquad p \le b .
\end{equation}
Applying a first-order Taylor expansion to $S_v$ around $\mathbf x$, then evaluating at $\bm 0$ gives:
\begin{equation}\label{eq.approx}
    S_v(\bm 0)  \approx S_v(\mathbf x) + \sum_{ p \le b} w^v_p(\mathbf x) \, (0 - x_p).
\end{equation} 
Subsequently, we obtain
\begin{equation}\label{eq.approx2}
    S_v(\mathbf x) \approx S_v(\bm 0) + \sum_{ p \le b} w^v_p(\mathbf x) \,  x_p.
\end{equation}
The quality of the approximation in equation~\eqref{eq.approx2} is examined in Section~\ref{ss.taylor}.
Let $E$ be an elliptic curve over $\mathbb{Q}$.
After setting $\mathbf x=(\widetilde{a_p}(E))_{p \le b}$, equation~\eqref{eq.approx2} suggests that the saliency scores $w_p^v(\mathbf x)$ can be viewed as somewhat
analogous to the weights in the Mestre--Nagao sum:
\begin{equation}\label{eq.MNsum}
\frac{1}{\log(b)}\sum_{p \le b}\frac{\log(p)}{p}\, a_p(E)=\frac{2}{\log(b)}\sum_{p \le b}\frac{\log(p)}{\sqrt{p}} \, \widetilde{a}_p(E), \ \quad b\in\mathbb{R}.\end{equation}
In Section \ref{ss.matched} we visualise averaged saliency scores.
By doing that, we are able to observe a visual comparison with Mestre--Nagao sums and interpret the CNN predictions.

\subsection{Accuracy of the linear approximation}\label{ss.taylor}
We evaluate the residual
\begin{equation}\label{eq.residual}
    R_v(\mathbf x) = S_v(\mathbf x) - S_v(\bm 0) - \sum_{p \le b} w_p^v(\mathbf x)\, x_p ,
\end{equation}
with $v = c(\mathbf x)$ the predicted class, on a fixed sample of $1000$ test curves, for every training epoch and every conductor interval considered in Section~\ref{ss.matched}.
We report two normalisations of the residual. 
The first normalisation is relative to the magnitude of the score change being approximated: $|R_v(\mathbf x)|/|S_v(\mathbf x)-S_v(\bm 0)|$. 
The second normalisation is relative to the decision margin of the classifier: $|R_v(\mathbf x)|/(S_{(1)}(\mathbf x)-S_{(2)}(\mathbf x))$, where $S_{(1)} \geq S_{(2)}$ denote the two largest class scores.

\begin{figure}[htbp]
    \centering
    \includegraphics[width=0.95\textwidth]{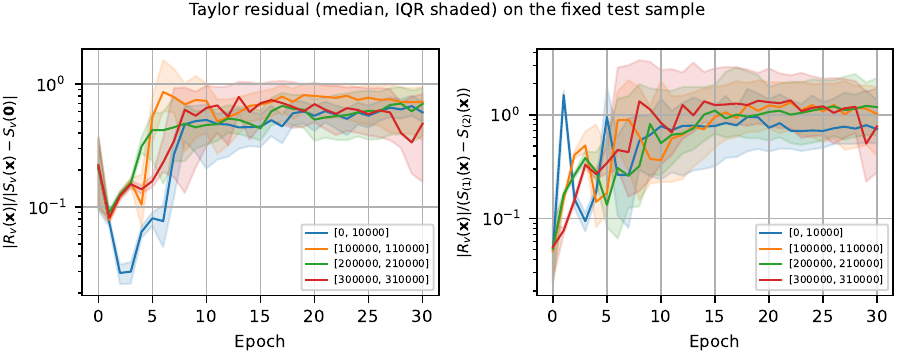}
    \caption{Median of the normalized Taylor residuals of equation~\eqref{eq.residual} over a fixed sample of $1000$ test curves, as a function of the training epoch, for the four conductor intervals, with interquartile range shaded. Left: residual relative to $|S_v(\mathbf x)-S_v(\bm 0)|$. Right: residual relative to $S_{(1)}(\mathbf x)-S_{(2)}(\mathbf x)$, that is, the top two score margin.}
    \label{fig:residual}
\end{figure}

The results are shown in Figure~\ref{fig:residual}.
Across the four conductor intervals and all training epochs, the median of the score-normalised residual ranges from roughly $0.3$ to $1.0$ 
and the margin-normalised residual is of comparable size.
In early epochs, this linear interpretation of the saliency scores is very accurate.
In the first 3 epochs, in all four conductor intervals, the median residual
$\lvert R_v(\mathbf{x})\rvert$ is at most $2.4\%$ of the difference of class
scores $\max_i S_i(\mathbf x) - \min_i S_i(\mathbf x)$
and the linearised scores $S_v(\bm 0)+\sum_{p \le b} w_p^v(\mathbf x)\,x_p$
attain their maximum at the network's predicted class for every curve in the sample.
As training continues, the pointwise approximation degrades.
This is visible in the initially small IQR bands in Figure~\ref{fig:residual} that grow with epochs.
One way to interpret the saliency curves is as coefficients of an approximately linear score during initial epochs, and as aggregate network scores in later epochs.

\subsection{Visualisation and interpretation}\label{ss.matched}
For an interval $[A,B]\subset\mathbb{R}$, we denote by \texttt{XECQ$^\circ$[A,B]} a fixed subset of $\numprint{3000}$ test curves and, for $v\in\{0,1,2,3\}$, we let \texttt{XECQ$^\circ$[A,B|v]} be the subset of $\texttt{XECQ$^\circ$[A,B]}$ with predicted rank $v$.
The same $\numprint{3000}$ curves are used at every epoch, so that the evolution of the saliency curves reflects the evolution of the network alone and not a change in the sample (see also Section~\ref{ss.control}).
For $E\in \texttt{XECQ$^\circ$[A,B]}$, we take $\mathbf x=(\widetilde{a_q}(E))_{q \le b}$ in equation~\eqref{eq.wp-1} and write $w_p^v (E) \seteq w_p^v (\mathbf x)$ to indicate that $\mathbf x=(\widetilde{a_q}(E))_{q \le b}$ is determined by $E$.

\begin{figure}
    \centering
    \includegraphics[width=0.5\linewidth]{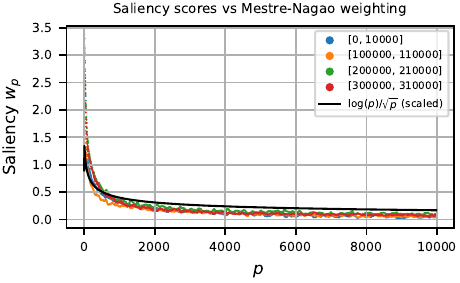}
    \caption{Comparison of the saliency scores $w_p$ with the Mestre--Nagao weighting. In this figure, the saliency scores are computed after $100$ epochs.}
    \label{fig:saliency_comparison}
\end{figure}

\begin{figure}
    \centering
    \includegraphics[width=0.5\linewidth]{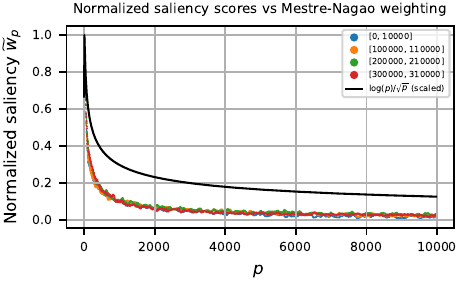}
    \caption{Comparison \rev{of} the normalized saliency scores $\widetilde{w}_p$ with the Mestre--Nagao weighting. In this figure, the normalised saliency scores are computed after $100$ epochs.}
    \label{fig:saliency_comparison_normalized}
\end{figure}

For each prime $p \le b$, define
\begin{equation}\label{eq.saliencyscore}
w_p \seteq \frac{1}{\#\texttt{XECQ$^\circ$[A,B]}}\sum_{v\in\{0,1,2,3,4\}}\sum_{E\in \texttt{XECQ$^\circ$[A,B|v]}} |w^v_p(E)|,
\end{equation}
and also consider the following normalisation of $w_p$:
\[\widetilde{w}_p=\frac{w_p}{\max\{w_q:q \le b, \ q \text{ prime}\}}\in[0,1].\]
In Figure~\ref{fig:saliency_comparison} (resp. Figure~\ref{fig:saliency_comparison_normalized}) we plot $w_p$ (resp. $\widetilde{w_p}$)  and $\log(p)/\sqrt{p}$, that is, the coefficient for $\widetilde{a_p}(E)$ in equation~\eqref{eq.MNsum}.

The averaged saliency score in equation~\eqref{eq.saliencyscore} sums over all possible ranks in order to compare with the Mestre--Nagao weightings, which are the same for every curve, irrespective of its analytic rank. 

In order to interpret the predictions made by the CNN, we need to understand the saliency of each class.
To that end, for $v\in\{0,1,2,3,4\}$, we define
\begin{equation}\label{eq.bigW}
W_p^v \seteq \frac{1}{\#\texttt{XECQ$^\circ$[A,B|v]}}\sum_{E\in \texttt{XECQ$^\circ$[A,B|v]}}w_p^v(E).
\end{equation}
It is conventional to take the absolute value of the saliency scores (as in equation~\eqref{eq.saliencyscore}), however the sign is perhaps interesting from an arithmetic perspective and so we not do that in equation~\eqref{eq.bigW}.
That being said, for each $v\in\{0,1,2,3,4\}$, we normalize $W_p^v$ as follows: 
\[\widetilde{W}_p^v=\frac{W_p^v}{\max\{|W_q^v|:q \le b, \ q \text{ prime}\}}\in[-1,1],\]
in which the denominator does utilise an absolute value.
Plotting \rev{$p/B$} against $\widetilde{W}_p^v$ after each epoch yields the so-called {\em saliency curve} for $v$.
These curves highlight which primes the network deems most important for its predictions.

\begin{figure}[htbp]
    \centering
    \includegraphics[width=1.0\textwidth]{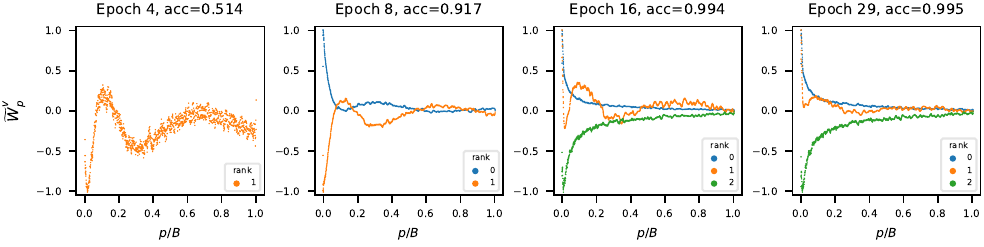}
    \caption{Saliency curves for \texttt{XECQ}$[0,10000]$, $p \le 10^4$, at four representative epochs, grouped by predicted rank. See Conventions~\ref{conv}.}
    \label{fig:sal_row_0}
\end{figure}

\begin{figure}[htbp]
    \centering
    \includegraphics[width=1.0\textwidth]{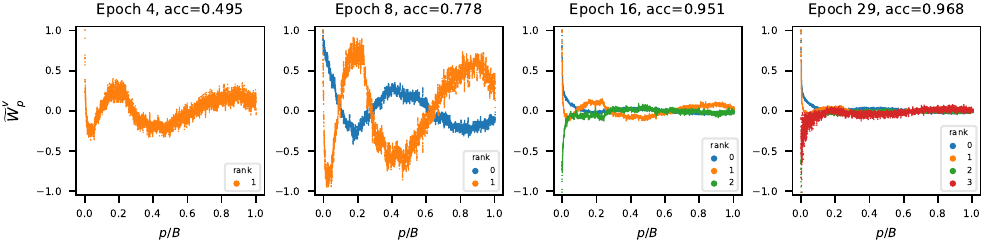}
    \caption{Saliency curves for \texttt{XECQ}$[100000,110000]$, $p \le 1.1\times 10^5$, grouped by predicted rank. See Conventions~\ref{conv}.}
    \label{fig:sal_row_1e5}
\end{figure}

\begin{figure}[htbp]
    \centering
    \includegraphics[width=1.0\textwidth]{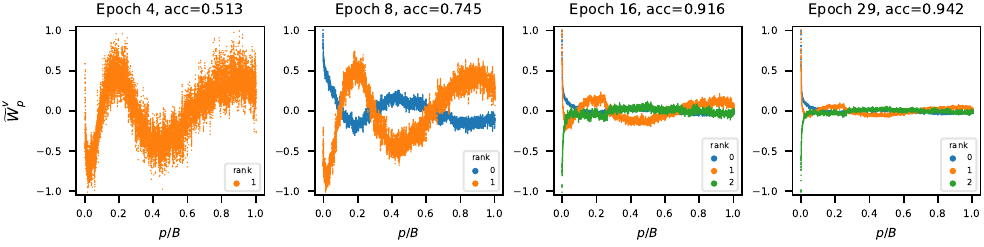}
    \caption{Saliency curves for \texttt{XECQ}$[200000,210000]$, $p \le 2.1\times 10^5$, grouped by predicted rank. See Conventions~\ref{conv}.}
    \label{fig:sal_row_2e5}
\end{figure}

\begin{figure}[htbp]
    \centering
    \includegraphics[width=1.0\textwidth]{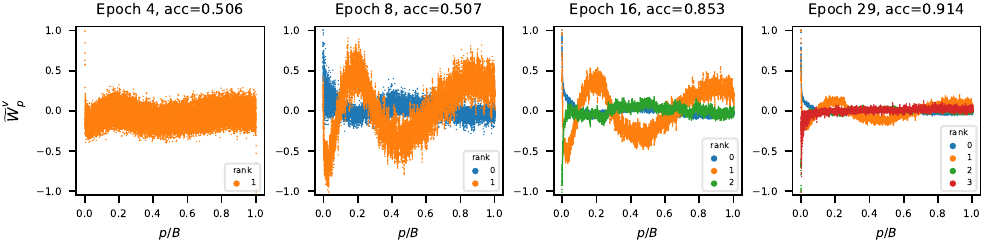}
    \caption{Saliency curves for \texttt{XECQ}$[300000,310000]$ at the matched murmuration scale $p \le 3.1\times 10^5$, grouped by predicted rank. Note the persistence of the murmuration oscillation through the final epoch. See Conventions~\ref{conv}.}
    \label{fig:sal_row_3e5}
\end{figure}

Figures \ref{fig:sal_row_0}, \ref{fig:sal_row_1e5}, \ref{fig:sal_row_2e5} and \ref{fig:sal_row_3e5} show the saliency curves at four representative epochs for the conductor intervals in Table~\ref{tab:rank_counts}.
The reader is referred to Conventions~\ref{conv} for the common features of the figures. 
The complete evolution over all $100$ epochs appears in \cite{code-repo}.
We observe that the saliency curves exhibit notable differences across conductor intervals $[A,B]$, ranks $v$, and training stages. 
 
For all four intervals, in the early epochs, the saliency curves for ranks $0$ and $1$ 
are in keeping with the murmuration phenomenon of \cite{HLOP}.
Remarkably, this behaviour emerges even though there is no visible murmuration in the average value of $\widetilde{a}_p$ \cite{drewLetter}, suggesting that the CNN is able to identify the hidden pattern in the dataset.
On the other hand, for $[0,10000]$, the oscillatory structure fades in later epochs. The saliency curve for rank $0$ (blue) gradually flattens and places increasing emphasis on smaller primes, resembling the Mestre--Nagao weightings in equation~\eqref{eq.MNsum}. 
Furthermore, a saliency curve for $v=2$ (green) appears, resembling the shape of the curve for $v=0$ but with the opposite sign.
These observations align with the pattern from \cite[Figure~1]{HLOP}, where rank $1$ is distinguished by the murmuration, and ranks $0$ and $2$ are separated mainly by the smaller primes.
For $[300000,310000]$, by contrast, the oscillation in the rank-$1$ saliency curve remains prominent throughout all epochs.


\begin{convention} \label{conv} 
Figures \ref{fig:sal_row_0}--\ref{fig:sal_row_true} (and the full grids in \cite{code-repo}) share the same conventions.
Each epoch consists of approximately $10$ steps (batches of $\numprint{3000}$ curves), the plots display the state of the network at step $0$ of the stated epoch, and the quoted accuracy is the test accuracy at that point.
Each color corresponds to a different predicted rank.
In the earliest epochs the network assigns a single predicted rank to every curve (which rank this is depends on the random initialization), so only one curve appears. 
Further ranks appear as training progresses.
The saliency scores are normalized so that their maximum absolute value is equal to $1$.
\end{convention}

\subsection{Grouping by true rank}\label{ss.control}
\begin{figure}[htbp]
    \centering
    \includegraphics[width=1.0\textwidth]{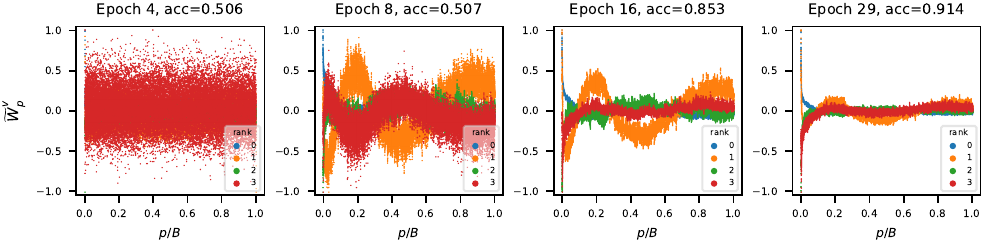}
    \caption{Saliency curves for \texttt{XECQ}$[300000,310000]$, averaged over curves grouped by their true analytic rank, on the same fixed test subset as Figure~\ref{fig:sal_row_3e5}.}
    \label{fig:sal_row_true}
\end{figure}

The sets \texttt{XECQ$^\circ$[A,B|v]} over which equation~\eqref{eq.bigW} averages are defined by the {\em predicted} rank, and these sets change as training progresses.
The evolution of the saliency curves could therefore in principle reflect the changing composition of the predicted classes rather than a change in the network gradients.
As a control, we recompute all saliency curves with the average in equation~\eqref{eq.bigW} taken over the curves whose true analytic rank is $v$ (the gradient is still taken at the predicted-class score), using the same fixed $\numprint{3000}$-curve subset at every epoch.
Figure~\ref{fig:sal_row_true} shows the result for the interval $[300000,310000]$. 
The corresponding figures for the other intervals appear in \cite{code-repo}.

\subsection{Early-phase saliency curve and murmurations}\label{ss.identification}

Combining \S\ref{ss.taylor} and \S\ref{ss.matched}, we can explain why murmuration-like patterns appear in the saliency curves.
In Section~\ref{ss.taylor}, we observed that, in early epochs, the network behaves like a linear score $S_v(\mathbf x) \approx S_v(\bm 0) + \sum_{p\le b} w_p^v\, x_p$.
For a linear score trained by gradient descent from small random initialization,
the first optimisation steps alter the weight vector along the difference of
class means: the gradient of the cross-entropy loss at near-zero weights is proportional to
$\sum_E \mathbf x_E\,(\mathbbm{1}_{[E \in v]} - P_v)$.
Stated differently, the early weight vectors tend to be approximately proportional
to the mean of $\mathbf{x_E}$ over the class $v$, minus the overall mean.
The $p$-th coordinate of the difference of class means is
\[
\mu_p \seteq \operatorname{avg}\{\widetilde a_p(E) : \operatorname{rank} E = 1\} - \operatorname{avg}\{\widetilde a_p(E) : \operatorname{rank} E = 0\},
\]
taken over a conductor interval.
This is precisely the empirical murmuration curve in the sense of~\cite{HLOP}.

\begin{figure}[htbp]
    \centering
    \includegraphics[width=0.85\textwidth]{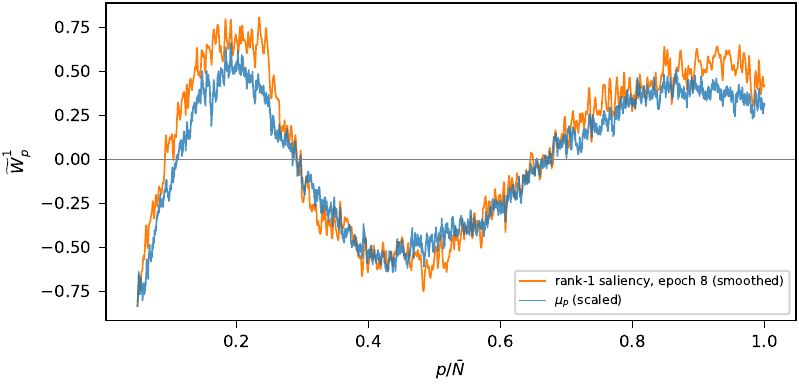}
    \captionof{figure}{The rank-$1$ true-rank saliency curve in early epochs and
    the empirical murmuration curve $(\mu_p)$ (scaled to match size) of the
    dataset, for the $b=B$ model on the conductor interval $[\numprint{100000},
    \numprint{110000}]$. (Points in the saliency curve are averaged over sliding windows of $25$ primes to reduce noise).} \label{fig:salvsmurm}
\end{figure}

Thus we should expect that the early-phase rank-$1$ saliency curve should
approximately coincide with $(\mu_p)_p$ up to scaling (and the rank-$0$ curve should approximately coincide with $-(\mu_p)_p$).
Figure~\ref{fig:salvsmurm} overlays the two curves for the interval $[\numprint{100000}, \numprint{110000}]$.
The match is clear enough that one could imagine predicting murmurations
from seeing these early epoch saliency curves (if one had the insight to examine these curves).
In hindsight we might have expected this, as the first principal directions (the
source of the initial murmuration observation in~\cite{HLOP}) and the major
class-mean differences are closely related.

\subsection{Check with synthetic data}
\begin{figure}[htbp]
    \centering
    \includegraphics[width=1.0\textwidth]{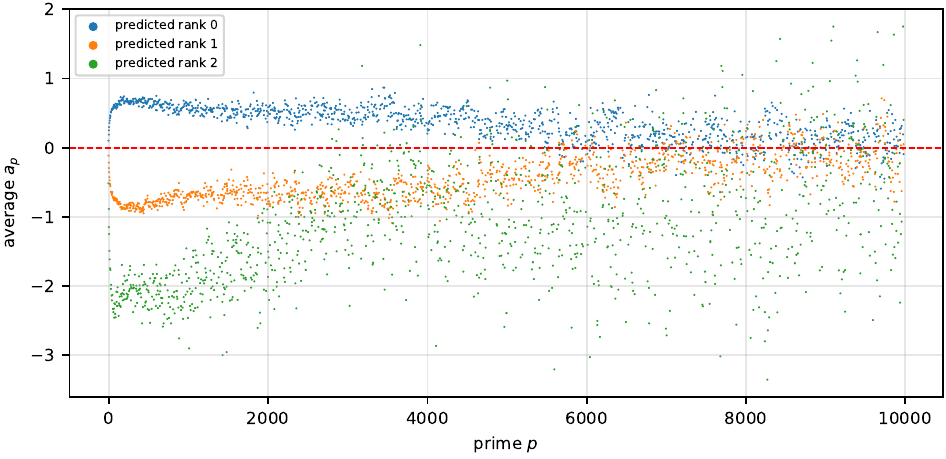}
    \captionof{figure}{Average values of $a_p$ (unnormalized) for the synthetic data, grouped by the rank predicted by the CNN. 
The network produces a clear separation between the predicted classes.}
    \label{fig:fake_ec} 
\end{figure}

To further investigate what CNN has learned, we conducted an experiment using synthetic data. Use of synthetic data was adopted in \cite{HLOa}, to which we refer the reader for the discussion of the Sato--Tate distribution of $(a_p)_p$. 
We generated a dataset of 500,000 sequences $(a_p)_p$ by drawing each coefficient independently from the Sato--Tate distribution, thereby producing data that mimics the statistical behaviour of elliptic
$L$-functions but is highly unlikely to arise from an elliptic curve in our dataset.  It is possible to identify each synthetic datapoint with some elliptic curve of large height and, therefore, most likely of much larger conductor. 
We then fed these synthetic sequences into the CNN trained on \texttt{XECQ} with conductor range $[0, 10000]$ and recorded the predicted rank for each sequence.
Finally, we grouped the synthetic curves by their predicted rank and computed the average unnormalized $a_p$ within each group.
The results are shown in Figure \ref{fig:fake_ec}.
Remarkably, the average $a_p$ for sequences predicted to be rank $0$, $1$ and $2$ exhibit a clear separation, which resembles the divergence observed in real data.
This indicates that the CNN associates certain bias patterns in the Frobenius traces with its rank predictions.

\clearpage

\printbibliography

\end{document}